\documentstyle[12pt]{article}
\begin{document}
\begin{titlepage}
\begin{flushright}
UQMATH-qe2-9702\\
q-alg/9703027
\end{flushright}
\vskip.3in

\begin{center}
{\Large \bf Super Yangian Double and its Central Extension}
\vskip.3in
{\large Yao-Zhong Zhang} \footnote{Queen Elizabeth II Fellow;
Email: yzz@maths.uq.edu.au}
\vskip.2in
{\em Department of Mathematics, University of Queensland, Brisbane,
Qld 4072, Australia}
\end{center}

\vskip 2cm
\begin{center}
{\bf Abstract}
\end{center}
We introduce super Yangian double $DY_\hbar[gl(m|n)]$ and its
central extension $\widehat{DY_\hbar[gl(m|n)]}$. We give their
defining relations in terms of current generators and obtain
Drinfeld comultiplication.


\end{titlepage}


\def\a{\alpha}
\def\b{\beta}
\def\d{\delta}
\def\e{\epsilon}
\def\g{\gamma}
\def\k{\kappa}
\def\l{\lambda}
\def\o{\omega}
\def\t{\theta}
\def\s{\sigma}
\def\D{\Delta}
\def\L{\Lambda}

\def\G{{\cal G}}
\def\C{{\bf C}}
\def\P{{\bf P}}

\def\dyglmn{\widehat{DY_\hbar[gl(m|n)]}}
\def\dygl11{\widehat{DY_\hbar[gl(1|1)]}}


\def\beq{\begin{equation}}
\def\eeq{\end{equation}}
\def\bea{\begin{eqnarray}}
\def\eea{\end{eqnarray}}
\def\ba{\begin{array}}
\def\ea{\end{array}}
\def\no{\nonumber}
\def\lt{\left}
\def\rt{\right}
\newcommand{\bq}{\begin{quote}}
\newcommand{\eq}{\end{quote}}

\newtheorem{Theorem}{Theorem}
\newtheorem{Definition}{Definition}
\newtheorem{Proposition}{Proposition}
\newtheorem{Lemma}[Theorem]{Lemma}
\newtheorem{Corollary}[Theorem]{Corollary}
\newcommand{\proof}[1]{{\bf Proof. }
        #1\begin{flushright}$\Box$\end{flushright}}

This paper concerns with Drinfeld current realization \cite{Dri88} of
super Yangian double $DY_\hbar[gl(m|n)]$ and its central extension
$\widehat{DY_\hbar[gl(m|n)]}$.

The Yangian double \cite{Kho96a} $DY_\hbar({\cal G})$ of a simple bosonic 
Lie algebra ${\cal G}$ is a quantum double of the
Yangian $Y_\hbar({\cal G})$ \cite{Dri88}. It is a deformation of the 
entire loop algebra and has important applications in
massive integrable models \cite{Ber93,Smi92}. 

The Yangian double with center (or central extension of the Yangian
double) $\widehat{DY_\hbar({\cal G})}$ for ${\cal G}=gl(n),~sl(n)$ 
were introduced in \cite{Kho96b,Ioh96} in terms of Drinfeld current
generators.

The philosophy behind this paper is to introduce
super Yangian double $DY_\hbar[gl(m|n)]$ and its central extension
$\dyglmn$. This is achieved by generalizing the
Reshetikhin and Semenov-Tian-Shansky
(RS) construction \cite{Res90} to the supersymmetric case. Using this
super RS construction and Gauss decomposition
\cite{Din93}, we obtain the defining relations for $\dyglmn$
in terms of super current generators. The computation in this paper
is parallel to our
recent work \cite{Zha97} on Drinfeld current
realization of quantum affine superalgebra $U_q[gl(m|n)^{(1)}]$
(see also \cite{Cai97a} for the special case of $m=n=1$), which in some
sense is a superization of work \cite{Fan97}.

The graded Yang-Baxter equation (YBE) with spectral-parameter dependence takes 
the form
\beq\label{rrr}
R_{12}({u-v})R_{13}(u)R_{23}(v)=R_{23}(v)R_{13}(u)R_{12}({u-v}),
\eeq
where $R(u)\in End(V\otimes V)$ with $V$ being graded vector space
and  obeys the weight conservation condition: 
$R(u)^{\a'\b'}_{\a\b}\neq 0$  only when
$[\a']+[\b']+[\a]+[\b]=0$ mod$2$.
The multiplication rule for the tensor product is defined for
homogeneous elements $a,~b,~c,~d$ of a quantum superalgebra by
\beq
(a\otimes b)(c\otimes d)=(-1)^{[b][c]}\,(ac\otimes bd)
\eeq
where $[a]\in {\bf Z}_2$ denotes the grading of the element $a$. 

Introduce the graded permutation operator $P$ on the tensor product
module $V\otimes V$ such that $P(v_\a\otimes v_\b)=(-1)^{[\a][\b]}
(v_\b\otimes v_\a)\,,~\forall v_\a, v_\b\in V$. 
In most cases R-matrix enjoys, among others, the following properties 
\bea
&& (i)~ P_{12}R_{12}(u)P_{12}=R_{21}(u),\label{pt-symmetry}\\
&& (ii)~ R_{12}(u)R_{21}(-u)=1.\label{unitarity}
\eea

The graded YBE, when written in matrix form, 
carries extra signs \cite{Bra94,Zha96},
\bea
&&R({u-v})_{\a\b}^{\a'\b'}
  R(u)_{\a'\g}^{\a''\g'}
  R(v)_{\b'\g'}^{\b''\g''}
  (-1)^{[\a][\b]+[\g][\a']+[\g'][\b']}\nonumber\\
&&~~~~~~=R(v)_{\b\g}^{\b'\g'}
  R(u)_{\a\g'}^{\a'\g''}
  R({u-v})_{\a'\b'}^{\a''\b''}
  (-1)^{[\b][\g]+[\g'][\a]+[\b'][\a']}.\label{ybe-def2'}
\eea

In \cite{Zha97} we generalized the RS construction \cite{Res90} 
to supersymmetric cases and obtained Drinfeld realization of the
quantum affine superalgebra $U_q[gl(m|n)^{(1)}]$.
Here we consider rational solution $R(u)$  to the graded YBE and give
a `rational' super RS algebra:

\begin{Definition}\label{rs}:  Rational super RS algebra is
generated by invertible $L^\pm(u)$, satisfying
\bea
R({u-v})L_1^\pm(u)L_2^\pm(v)&=&L_2^\pm(v)L_1^\pm(u)R({u-v}),\no\\
R({u_+-v_-})L_1^+(u)L_2^-(v)&=&L_2^-(v)L_1^+(u)R({u_--
         v_+}),\label{super-rs}
\eea
where $L_1^\pm(u)=L^\pm(u)\otimes 1$, $L_2^\pm(u)=1\otimes L^\pm(u)$
and $u_\pm=u\pm \frac{1}{2}\hbar c$. For the first formula of
(\ref{super-rs}), the expansion direction of $R({u-v})$ can be
chosen in $u\over v$ or $v\over u$, but for the second formula, the
expansion direction must only be in $u\over v$.
\end{Definition}

In matrix form, (\ref{super-rs}) carries extra signs due to the graded
multiplication rule of tensor products:
\bea
&&R({u-v})_{\a\b}^{\a"\b"}L^\pm(u)_{\a"}^{\a'}L^\pm(v)_{\b"}
      ^{\b'}\,(-1)^{[\a']([\b']+[\b"])}\no\\
&&~~~~~~~~~~~~~~~=L^\pm(v)_\b^{\b"}L^\pm(u)_\a^{\a"}R({u-v})
      _{\a"\b"}^{\a'\b'}\,(-1)^{[\a]([\b]+[\b"])},\no\\
&&R({u_+-v_-})_{\a\b}^{\a"\b"}L^+(u)_{\a"}^{\a'}L^-(v)_{\b"}
      ^{\b'}\,(-1)^{[\a']([\b']+[\b"])}\no\\
&&~~~~~~~~~~~~~~~=L^-(v)_\b^{\b"}L^+(u)_\a^{\a"}R({u_--v_+})
      _{\a"\b"}^{\a'\b'}\,(-1)^{[\a]([\b]+[\b"])}.\label{rll-component}
\eea
Introduce matrix $\t$:
\beq
\t^{\a'\b'}_{\a\b}=(-1)^{[\a][\b]}\d_{\a\a'}\d_{\b\b'}.\label{theta}
\eeq
With the help of this matrix $\t$, one can cast (\ref{rll-component}) into the
usual matrix equation,
\bea
R({u-v})L_1^\pm(u)\t L_2^\pm(v)\t&=&\t L_2^\pm(v)\t L_1^\pm(u)R({u-
         v}),\no\\
R({u_+-v_-})L_1^+(u)\t L_2^-(v)\t &=&\t L_2^-(v)\t L_1^+(u)R({u_--
         v_+}).\label{RLL-LLR1}
\eea
Now the multiplications in ({\ref{RLL-LLR1}) are simply the usual matrix
multiplications.

We will take $R({u})\in End(V\otimes V)$ to be the 
Yang's rational R-matrix associated with superalgebra $gl(m|n)$,
\beq
R(u)=\frac{1}{u+2\hbar}(uI+2\hbar P),
\eeq
where $V$ is a $(m+n)$-dimensional graded vector space. Let
basis vectors $\{v_1,\;v_2,\cdots,\;v_m\}$ be even and $\{v_{m+1},\;v_{m+2},
\cdots,\; v_{m+n}\}$ be odd. Then the R-matrix has the following matrix 
elements:
\bea
R({u})_{\a\b}^{\a'\b'}&=&(-1)^{[\a][\b]}\tilde{R}({u})_{\a\b}
  ^{\a'\b'},\no\\
\tilde{R}({u})&=&\sum^m_{i=1}E^i_i\otimes E^i_i
  +\sum_{i=m+1}^{m+n}\frac{2\hbar-u}{u+2\hbar}E^i_i\otimes E^i_i
  +\frac{u}{u+2\hbar}\sum_{i\neq j}(-1)^{[i][j]}E^i_i\otimes E^j_j\no\\
& &\sum_{i<j} \frac{2\hbar}{u+2\hbar}E^j_i\otimes E^i_j
  +\sum_{i>j}\frac{2\hbar}{u+2\hbar}E^j_i\otimes E^i_j.\label{r-matrix}
\eea
It is easy to check that the R-matrix $R(u)$ satisfies
(\ref{pt-symmetry}) and (\ref{unitarity}).
We will construct central extended super Yangian double $\dyglmn$.  

\begin{Theorem}\label{df}: $L^\pm(u)$ has the following
Gauss decomposition
\bea
L^\pm(u)&=&\left (
\begin{array}{cccc}
1 & \cdots & {} & 0\\
e^\pm_{2,1}(u) & \ddots & {} & {}\\
e^\pm_{3,1}(u) & {}     & {} & \vdots\\
\vdots &  {} & {} & {}\\
e^\pm_{m+n,1}(u) & \cdots & e^\pm_{m+n,m+n-1}(u) & 1
\end{array}
\right )
\left (
\begin{array}{ccc}
k^\pm_1(u) & \cdots & 0\\
\vdots & \ddots & \vdots\\
0 & \cdots & k^\pm_{m+n}(u)
\end{array}
\right )\no\\
& & \times \left (
\begin{array}{ccccc}
1 & f^\pm_{1,2}(u) & f^\pm_{1,3}(u) & \cdots & f^\pm_{1,m+n}(u)\\
\vdots & \ddots & \cdots & {} & \vdots\\
{} & {} & {} & {} & f^\pm_{m+n-1,m+n}(u)\\
0 & {} & {} & {} & 1
\end{array}
\right ),
\eea
where $e^\pm_{i,j}(u),~f^\pm_{j,i}(u)$ and $k^\pm_i(u) ~(i>j)$ are 
elements in the rational super RS algebra and $k^\pm_i(u)$ are invertible. 
Let
\bea
X^-_i(u)&=&f^+_{i,i+1}(u_+)-f^-_{i,i+1}(u_-),\no\\
X^+_i(u)&=&e^+_{i+1,i}(u_-)-e^-_{i+1,i}(u_+),
\eea
where $u_\pm=u\pm\frac{1}{2}\hbar c$, then $c,\;
X^\pm_i(u),\;k^\pm_j(u),~i=1,2,\cdots,m+n-1,\;j=1,2,\cdots,m+n$ give
the defining relations of central extended super Yangian double  $\dyglmn$
which, when $c=0$ reduce to those of super Yangian double
$DY_\hbar[gl(m|n)]$.
\end{Theorem}

The Gauss decomposition implies that the elements $e^\pm_{i,j}(u),\;
f^\pm_{j,i}(u) ~(i>j)$ and $k^\pm_i(u)$ are uniquely determined by
$L^\pm(u)$. In the following we will denote $f^\pm_{i,i+1}(u),\;
e^\pm_{i+1,i}(u)$ as $f^\pm_i(u),\;e^\pm_i(u)$, respectively.

The following matrix equations can be deduced from (\ref{RLL-LLR1}):
\bea
R_{21}({u-v})\t L^\pm_2(u)\t L_1^\pm(v)&=&
  L^\pm_1(v)\t L^\pm_2(u)\t R_{21}({u-v}),\label{RLL-LLR2}\\
R_{21}({u_--v_+})\t L^-_2(u)\t L_1^+(v)&=&
  L^+_1(v)\t L^-_2(u)\t R_{21}({u_+-v_-}),\label{RLL-LLR3}\\
\t L^\pm_2(u)^{-1}\t L^\pm_1(v)^{-1}R_{21}({u-v})&=&
  R_{21}({u-v})L^\pm_1(v)^{-1}\t L^\pm_2(u)^{-1}
  \t,\label{RLL-LLR4}\\
\t L^+_2(u)^{-1}\t L^-_1(v)^{-1}R_{21}({u_+-v_-})&=&
  R_{21}({u_--v_+})L^-_1(v)^{-1}\t L^+_2(u)^{-1}
  \t,\label{RLL-LLR5}\\
L^\pm_1(v)^{-1}R_{21}({u-v})\t L^\pm_2(u)\t&=&
  \t L^\pm_2(u)\t R_{21}({u-
  v})L^\pm_1(v)^{-1},\label{RLL-LLR6}\\
L^-_1(v)^{-1}R_{21}({u_+-v_-})\t L^+_2(u)\t&=&
  \t L^+_2(u)\t R_{21}({u_--
  v_+})L^-_1(v)^{-1},\label{RLL-LLR7}\\
L^+_1(v)^{-1}R_{21}({u_--v_+})\t L^-_2(u)\t&=&
  \t L^-_2(u)\t R_{21}({u_+-
  v_-})L^+_1(v)^{-1},\label{RLL-LLR8}
\eea
where $R_{21}({u-v})=R({v-u})^{-1}$. As in (\ref{RLL-LLR1}), the
multiplications in (\ref{RLL-LLR2} -- \ref{RLL-LLR8}) are usual matrix
multiplications.

Using (\ref{theta}), (\ref{r-matrix}), (\ref{RLL-LLR1},
\ref{RLL-LLR2} -- \ref{RLL-LLR8}) and theorem \ref{df}, and by 
parallel calculations as to \cite{Zha97}, we obtain

\begin{Definition}\label{general}: $\dyglmn$ is an
associative algebra over the ring of formal power series in the
variable $\hbar$ and with Drinfeld current generators: $X^\pm_i(u),~
k^\pm_j(u),~i=1,2,\cdots,m+n-1,~j=1,2,\cdots,m+n$ and a central element
$c$. $k^\pm_i(u)$ are invertible.
The grading of the generators are: $[X^\pm_m(u)]=1$ and zero otherwise.
When $c=0$, $\dyglmn$ reduces to $DY_\hbar[gl(m|n)]$.
The defining relations are given by
\bea
k^\pm_i(u)k^\pm_j(v)&=&k^\pm_j(v)k^\pm_i(u),~~i\neq j\no\\
k^+_i(u)k^-_i(v)&=&k^-_i(v)k^+_i(u),~~i\leq m,\no\\
{{u_+-v_--2\hbar}\over{u_+-v_-+2\hbar}}k^+_i(u)k^-_i(v)&=&
  {{u_--v_+-2\hbar}\over{u_--v_++2\hbar}}k^-_i(v)k^+_i(u),~~m<i\leq
  m+n,\no\\
{{u_\pm-v_\mp}\over{u_\pm -v_\mp +2\hbar}}k^\mp_i(v)^{-1}k^\pm_j(u)&=&
  {{u_\mp-v_\pm}\over{u_\mp -v_\pm
  +2\hbar}}k^\pm_j(u)k^\mp_i(v)^{-1},~~i>j,\no\\
k^\pm_j(u)^{-1}X^-_i(v)k^\pm_j(u)&=&X^-_i(v),~~j-i\leq -1,\no\\
k^\pm_j(u)^{-1}X^+_i(v)k^\pm_j(u)&=&X^+_i(v),~~j-i\leq -1,~~{\rm or}\no\\
k^\pm_j(u)^{-1}X^-_i(v)k^\pm_j(u)&=&X^-_i(v),~~j-i\geq 2,\no\\
k^\pm_j(u)^{-1}X^+_i(v)k^\pm_j(u)&=&X^+_i(v),~~j-i\geq 2,\no\\
k^\pm_i(u)^{-1}X^-_i(v)k^\pm_i(u)&=&
  \frac{u_\mp -v+2\hbar}{u_\mp-v}X^-_i(v),~~i<m,\no\\
k^\pm_i(u)^{-1}X^-_i(v)k^\pm_i(u)&=&
  \frac{u_\mp -v-2\hbar}{u_\mp-v}X^-_i(v),~~m<i\leq m+n-1,\no\\
k^\pm_{i+1}(u)^{-1}X^-_i(v)k^\pm_{i+1}(u)&=&
  \frac{u_\mp-v-2\hbar}{u_\mp-v}X^-_i(v),~~i<m,\no\\
k^\pm_{i+1}(u)^{-1}X^-_i(v)k^\pm_{i+1}(u)&=&
  \frac{u_\mp -v+2\hbar}{u_\mp-v}X^-_i(v),~~m<i\leq m+n-1,\no\\
k^\pm_i(u)X^+_i(v)k^\pm_i(u)^{-1}&=&
  \frac{u_\pm -v+2\hbar}{u_\pm-v}X^+_i(v),~~i<m,\no\\
k^\pm_i(u)X^+_i(v)k^\pm_i(u)^{-1}&=&
  \frac{u_\pm -v-2\hbar}{u_\pm-v}X^+_i(v),~~m<i\leq m+n-1,\no\\
k^\pm_{i+1}(u)X^+_i(v)k^\pm_{i+1}(u)^{-1}&=&
  \frac{u_\pm -v-2\hbar}{u_\pm-v}X^+_i(v),~~i<m,\no\\
k^\pm_{i+1}(u)X^+_i(v)k^\pm_{i+1}(u)^{-1}&=&
  \frac{u_\pm -v+2\hbar}{u_\pm-v}X^+_i(v),~~m<i\leq m+n-1,\no\\
k^\pm_i(u)^{-1}X^-_m(v)k^\pm_i(u)&=&
  \frac{u_\mp -v+2\hbar}{u_\mp-v}X^-_m(v),~~i=m,\;m+1,\no\\
k^\pm_i(u)X^+_m(v)k^\pm_i(u)^{-1}&=&
  \frac{u_\pm -v+2\hbar}{u_\pm-v}X^+_m(v),~~i=m,\;m+1,\no\\
(u-v\mp 2\hbar)X^\mp_i(u)X^\mp_i(v)&=&(u-v\pm 2\hbar)
  X^\mp_i(v)X^\mp_i(u),~~i<m,\no\\
(u-v\pm 2\hbar)X^\mp_i(u)X^\mp_i(v)&=&(u-v\mp 2\hbar)
  X^\mp_i(v)X^\mp_i(u),~~m<i\leq m+n-1,\no\\
\{X^\pm_m(u),X^\pm_m(v)\}&=&0,\no\\
(u-v)X^+_i(u)X^+_{i+1}(v)&=&(u-v+2\hbar)
  X^+_{i+1}(v)X^+_i(u),~~i<m,\no\\
(u-v)X^+_i(u)X^+_{i+1}(v)&=&(u-v-2\hbar)
  X^+_{i+1}(v)X^+_i(u),~~m\leq i\leq m+n-1,\no\\
(u-v+2\hbar)X^-_i(u)X^-_{i+1}(v)&=&(u-v)
  X^-_{i+1}(v)X^-_i(u),~~i<m,\no\\
(u-v-2\hbar)X^-_i(u)X^-_{i+1}(v)&=&(u-v)
  X^-_{i+1}(v)X^-_i(u),~~m\leq i\leq m+n-1,\no\\
{[X^+_i(u),X^-_j(v)]}&=&-2\hbar\d_{ij}\lt(
   \d(u_--v_+)k^+_{i+1}(v_+)
   k^+_i(v_+)^{-1}\rt.\no\\
& &~~~~~\lt. -\d(u_+-v_-)k^-_{i+1}(u_+)k^-_i(u_+)^{-1}
  \rt),~~i,j\neq m,\no\\
\{X^+_m(u),X^-_m(v)\}&=&2\hbar\lt(
   \d(u_--v_+)k^+_{m+1}(v_+)
   k^+_m(v_+)^{-1}\rt.\no\\
& &~~~~\lt. -\d(u_+-v_-)k^-_{m+1}(u_+)k^-_m(u_+)^{-1}\rt),
  \label{glmn-rs}
\eea
where $[X,Y]\equiv XY-YX$ stands for a commutator and $\{X,Y\}\equiv XY+YX$
for an anti-commutator and
\beq
\d(u-v)=\sum_{k\in {\bf Z}}\,u^kv^{-k-1}
\eeq
is a formal series, together with
the following Serre and extra Serre \cite{Sch92,Yam96} relations:
\bea
&&\{X^\pm_i(u_1)X^\pm_i(u_2)X^\pm_{i+1}(v)-2X^\pm_i(u_1)
  X^\pm_{i+1}(v)X^\pm_i(u_2)\no\\
&&~~~~~~~~~~  +X^\pm_{i+1}(v)X^\pm_i(u_1)X^\pm
  _i(u_2)\}+\{u_1\leftrightarrow u_2\}=0,~~i\neq m,
  \label{serre1}\\
&&\{X^\pm_{i+1}(u_1)X^\pm_{i+1}(u_2)X^\pm_{i}(v)-2X^\pm_{i+1}(u_1)
  X^\pm_{i}(v)X^\pm_{i+1}(u_2)\no\\
&&~~~~~~~~~~  +X^\pm_{i}(v)X^\pm_{i+1}(u_1)X^\pm
  _{i+1}(u_2)\}+\{u_1\leftrightarrow u_2\}=0,~~i\neq m-1,
  \label{serre2}\\
&&\{(u_1-u_2\mp 2\hbar)[X^\pm_m(u_1)X^\pm_m(u_2)X^\pm_{m-1}(v)-
  2X^\pm_m(u_1)X^\pm_{m-1}(v)X^\pm_m(u_2)\no\\
&&~~~~~~~~~~  +X^\pm_{m-1}(v)X^\pm_m(u_1)X^\pm
    _m(u_2)]\}+\{u_1\leftrightarrow u_2\}=0,
    \label{serre3}\\
&&\{(u_2-u_1\mp 2\hbar)[X^\pm_m(u_1)X^\pm_m(u_2)X^\pm_{m+1}(v)-
  2X^\pm_m(u_1)X^\pm_{m+1}(v)X^\pm_m(u_2)\no\\
&&~~~~~~~~~~  +X^\pm_{m+1}(v)X^\pm_m(u_1)X^\pm
     _m(u_2)]\}+\{u_1\leftrightarrow u_2\}=0,
     \label{serre4}\\
&&\{(u_1-u_2\mp 2\hbar)[X^\pm_m(u_1)X^\pm_m(u_2)X^\pm_{m-1}(v_1)
  X^\pm_{m+1}(v_2)\no\\
&&~~~~~~~~~~  -2X^\pm_m(u_1)X^\pm_{m-1}(v_1)X^\pm_m(u_2)
  X^\pm_{m+1}(v_2)]\no\\
&&~~~~~~~~~~\mp 4\hbar X^\pm_{m-1}(v_1)X^\pm_m(u_1)X^\pm
            _m(u_2)X^\pm_{m+1}(v_2)\no\\
&&~~~~~~~~~~+(u_2-u_1\mp 2\hbar)[-2X^\pm_{m-1}(v_1)
  X^\pm_m(u_1)X^\pm_{m+1}(v_2)X^\pm_m(u_2)\no\\
&&~~~~~~~~~~+X^\pm_{m-1}(v_1)X^\pm_{m+1}(v_2)X^\pm_m(u_1)
  X^\pm_m(u_2)]\}+\{u_1\leftrightarrow u_2\}=0.
  \label{extra-serre}
\eea
\end{Definition}

\noindent{\bf Remark}: For the special case of $m=n=1$, we have
\bea
k^\pm_i(u)k^\pm_j(v)&=&k^\pm_j(v)k^\pm_i(u),~~i,\;j=1,2,\no\\
k^+_1(u)k^-_1(v)&=&k^-_1(v)k^+_1(u),\no\\
{{u_+-v_--2\hbar}\over{u_+-v_-+2\hbar}}k^+_2(u)k^-_2(v)&=&
  {{u_--v_+-2\hbar}\over{u_--v_++2\hbar}}k^-_2(v)k^+_2(u),\no\\
{{u_\pm-v_\mp}\over{u_\pm-v_\mp+2\hbar}}k^\mp_2(v)^{-1}k^\pm_1(u)&=&
  {{u_\mp-v_\pm}\over{u_\mp-v_\pm
  +2\hbar}}k^\pm_1(u)k^\mp_2(v)^{-1},\no\\
k^\pm_i(u)^{-1}X^-_1(v)k^\pm_i(u)&=&
  \frac{u_\mp-v+2\hbar}{u_\mp-v}X^-_1(v),\no\\
k^\pm_i(u)X^+_1(v)k^\pm_i(u)^{-1}&=&
  \frac{u_\pm -v+2\hbar}{u_\pm-v}X^+_1(v),\no\\
\{X^\pm_1(u),X^\pm_1(v)\}&=&0,\no\\
\{X^+_1(u),X^-_1(v)\}&=&2\hbar\lt(\d(u_--v_+)k^+_2(v_+)
   k^+_1(v_+)^{-1}\rt.\no\\
& &~~~~~\lt. -\d(u_+-v_-)k^-_2(u_+)k^-_1(u_+)^{-1}
   \rt).\label{gl11-rs}
\eea
This is the defining relations of $\widehat{DY_\hbar[gl(1|1)]}$, which, 
when $c=0$, reduce to those of  centraless super Yangian double 
$DY_\hbar[gl(1|1)]$ obtained in \cite{Cai97b}.

\begin{Theorem}\label{hopf-glmn}: The algebra $\dyglmn$
given by definition \ref{general} has a Hopf algebra structure, which
is given by the following formulae.\\
{\bf Coproduct} $\D$
\bea
&&\D(c)=c\otimes 1+1\otimes c,\no\\
&&\D(k^\pm_j(u))=k^\pm_j(u\pm\frac{1}{2}\hbar c_2)\otimes k^\pm_j
   (u\mp\frac{1}{2}\hbar c_1),~~j=1,2,\cdots,m+n,\no\\ 
&&\D(X^+_i(u))=X^+_i(u)\otimes 1+\psi_i(u+\frac{1}{2}\hbar c_1)\otimes
   X^+_i(u+\hbar c_1),\no\\
&&\D(X^-_i(u))=1\otimes X^-_i(u)+ X^-_i(u+\hbar c_2)\otimes
  \phi_i(u+\frac{1}{2}\hbar c_2),~~i=1,2,\cdots,m+n-1,\no\\
  \label{coproduct-glmn}
\eea
where $c_1=c\otimes 1,~c_2=1\otimes c,~\psi_i(u)=k^-_{i+1}(u)k^-_i(u)^{-1}$
and $\phi_i(u)=k^+_{i+1}(u)k^+_i(u)^{-1}$.\\
{\bf Counit} $\e$
\beq
\e(c)=0,~~~~\e(k^\pm_j(u))=1,~~~~\e(X^\pm_i(u))=0.\label{counit-glmn}
\eeq
{\bf Antipode} $S$
\bea
&&S(c)=-c,~~~~S(k^\pm_j(u))=k^\pm_j(u)^{-1},\no\\
&&S(X^+_i(u))=-\psi_i(u-\frac{1}{2}\hbar c)^{-1}X^+_i(u-\hbar c),\no\\
&&S(X^-_i(u))=-X^-_i(u-\hbar c)\phi_i(u-\frac{1}{2}\hbar c)^{-1}.
    \label{antipode-glmn}
\eea
\end{Theorem}


\vskip.3in
\noindent {\bf Acknowledgements.} This work is supported by
Australian Research Council, and in part by University of Queensland
New Staff Research Grant and External Support Enabling Grant.



\vskip.3in
\begin{thebibliography}{99}
\bibitem{Dri88} V.G. Drinfeld, Sov. Math. Dokl. {\bf 36} (1988) 212.
\bibitem{Kho96a} S. Khoroshkin, V. Tolstoy, Lett. Math. Phys. {\bf 36}
   (1996) 319.
\bibitem{Ber93} D. Bernard, A. LeClair, Nucl. Phys. {\bf B399} (1993)
   709.
\bibitem{Smi92} F.A. Smirnov, Int. J. Mod. Phys. {\bf A7 suppl. 1B}
   (1992) 813.
\bibitem{Kho96b} S. Khoroshkin, {\em Central extension of the Yangian 
   double}, preprint q-alg/9602031.
\bibitem{Ioh96} K. Iohara, M. Kohno, Lett. Math. Phys. {\bf 37}
   (1996) 319.
\bibitem{Res90} N.Yu. Reshetikhin, M.A. Semenov-Tian-Shansky,
   Lett. Math. Phys. {\bf 19} (1990) 133.
\bibitem{Din93} J. Ding, I.B. Frenkel, Commun. Math. Phys. {\bf 155}
   (1993) 277.
\bibitem{Zha97} Y.-Z. Zhang, {\em Comments on Drinfeld realization of
   quantum affine superalgebra $U_q[gl(m|n)^{(1)}]$ and its Hopf
   algebra structure}, preprint q-alg/9703020.
\bibitem{Cai97a} J.-F. Cai, S.K. Wang, K. Wu, W.-Z. Zhao, {\em Drinfeld
   realization of quantum affine superalgebra $U_q[\widehat{gl(1|1)}]$},
   preprint q-alg/9703022.
\bibitem{Fan97} H. Fan, B.-Y. Hou, K.-J. Shi, J. Math. Phys. 
   {\bf 38} (1997) 411.
\bibitem{Zha96} Y.-Z. Zhang, {\em On the graded Yang-Baxter and
    reflection equations}, Commun. Theor. Phys. (1996). 
\bibitem {Bra94}  A.J. Bracken, G.W. Delius, M.D. Gould, Y.-Z. Zhang,  
    J. Phys. {\bf A:} Math. Gen. {\bf 27} (1994) 6551.
\bibitem{Sch92} M. Scheunert, Lett. Math. Phys. {\bf 24} (1992) 173.
\bibitem{Yam96} H. Yamane, {\em On defining relations of affine Lie
    superalgebras and their quantized universal enveloping
    superalgebras}, preprint q-alg/9603015.
\bibitem{Cai97b} J.-F. Cai, G.-X. Ju, K. Wu, S.-K. Wang, {\em Super
    Yangian double $DY(gl(1|1))$ and its Gauss decomposition}, preprint
    q-alg/9701011.
\end{thebibliography}
\end{document}